\documentclass[12pt]{article}
\usepackage{subeqn}
\usepackage{graphicx}
\usepackage{ifpdf}
\ifpdf  \fi
\usepackage{amsfonts,amssymb}
\usepackage{mathptmx,helvet,courier,makeidx,multicol,footmisc}
\usepackage[numbers]{natbib}
\bibpunct{(}{)}{;}{a}{,}{,}
\usepackage[bookmarksnumbered=true, pdfauthor={Wen-Long Jin}]{hyperref}

\newcommand{\commentout}[1]{}

\newcommand{\ba}{\begin{array}}
        \newcommand{\ea}{\end{array}}
\newcommand{\bc}{\begin{center}}
        \newcommand{\ec}{\end{center}}
\newcommand{\bdm}{\begin{displaymath}}
        \newcommand{\edm}{\end{displaymath}}
\newcommand{\bds} {\begin{description}}
        \newcommand{\eds} {\end{description}}%17Apr01
\newcommand{\ben}{\begin{enumerate}}
        \newcommand{\een}{\end{enumerate}}
\newcommand{\beq}{\begin{equation}}
        \newcommand{\eeq}{\end{equation}}
\newcommand{\bfg} {\begin{figure}[h]}
        \newcommand{\efg} {\end{figure}}%Nov 5,99
\newcommand{\bi} {\begin {itemize}}
        \newcommand{\ei} {\end {itemize}}
\newcommand{\bqn}{\begin{eqnarray}}
        \newcommand{\eqn}{\end{eqnarray}}
\newcommand{\bqs}{\begin{eqnarray*}}
        \newcommand{\eqs}{\end{eqnarray*}}
\newcommand{\bsl} {\begin{slide}[8.8in,6.7in]}
        \newcommand{\esl} {\end{slide}}
\newcommand{\bsq}{\begin{subequations}}
        \newcommand{\esq}{\end{subequations}}       
\newcommand{\bss} {\begin{slide*}[9.3in,6.7in]}
        \newcommand{\ess} {\end{slide*}}
\newcommand{\btb} {\begin {table}}
        \newcommand{\etb} {\end {table}}%Nov 10,99
\newcommand{\bvb}{\begin{verbatim}}
        \newcommand{\evb}{\end{verbatim}}

\newcommand{\m}{\mbox}
\newcommand {\der}[2] {{\frac {\m {d} {#1}} {\m{d} {#2}}}}

\newcommand {\pd}[2] {{\frac {\partial {#1}} {\partial {#2}}}}

\newcommand{\cas}[1]{{{\left \{ \ba #1 \ea \right. }}}

\newcommand{\reff}[1] {{{Figure \ref {#1}}}}
\newcommand{\refe}[1] {{(\ref {#1})}}%Nov 5

\def\pmb#1{\setbox0=\hbox{$#1$}%
   \kern-.025em\copy0\kern-\wd0
   \kern.05em\copy0\kern-\wd0
   \kern-.025em\raise.0433em\box0 }

\def\eop{{\hfill $\blacksquare$}}%17Apr01

\newtheorem{theorem}{Theorem}[section]%17Apr01
\newtheorem{definition}[theorem]{Definition}%17Apr01

\newtheorem{lemma}[theorem]{Lemma}%17Apr01
\def\dx     {{\Delta x}}
\def\dt     {{\Delta t}}

\usepackage[draft]{changes}
\definechangesauthor[name={Wenlong Jin},color=red]{WJ}

\usepackage{multirow}
\usepackage{subcaption}
\usepackage{float}
\usepackage[colorinlistoftodos]{todonotes}

\begin{document}
	\title{Stable day-to-day dynamics for departure time choice} %20170918
\author{Wen-Long Jin \footnote{Department of Civil and Environmental Engineering, California Institute for Telecommunications and Information Technology, Institute of Transportation Studies, 4000 Anteater Instruction and Research Bldg, University of California, Irvine, CA 92697-3600. Tel: 949-824-1672. Fax: 949-824-8385. Email: wjin@uci.edu. Corresponding author}}
\maketitle
\begin{abstract}
All existing day-to-day dynamics of departure time choice at a single bottleneck are unstable, and this has led to doubt regarding the existence of a stable user equilibrium in the real world. However, empirical observations and our personal driving experience suggest stable stationary congestion patterns during a peak period.
In this paper we attempt to reconcile the discrepancy by presenting a stable day-to-day dynamical system for drivers' departure time choice at a single bottleneck. 

We first define within-day traffic dynamics with the point queue model, costs,  the departure time user equilibrium (DTUE), and the arrival time user equilibrium (ATUE).
We then identify three behavioral principles: (i) Drivers choose their departure and arrival times in a backward fashion (backward choice principle); (ii) After choosing the arrival times, they update their departure times to balance the total costs (cost balancing principle); (iii) They choose their arrival times to reduce their scheduling costs or gain their scheduling payoffs (scheduling cost reducing or scheduling payoff gaining principle).
In this sense, drivers' departure and arrival time choices are driven by their scheduling payoff choice. With a single tube or imaginary road model, we convert the nonlocal day-to-day arrival time shifting problem to a local scheduling payoff shifting problem. After introducing a new variable for the imaginary density, we apply the Lighthill-Whitham-Richards (LWR) model to describe the day-to-day dynamics of scheduling payoff choice and present splitting and cost balancing schemes to determine arrival and departure flow-rates accordingly. We also develop the corresponding discrete models for numerical solutions. We theoretically prove that the day-to-day stationary state of the LWR model leads to the scheduling payoff user equilibrium (SPUE), which is equivalent to both DTUE and ATUE and stable. We use one numerical example to demonstrate the effectiveness and stability of the new day-to-day dynamical model. This study is the first step for understanding stable day-to-day dynamics for departure time choice, and many follow-up studies are possible and warranted.

\end{abstract}
 
\section{Introduction}
Congestion at critical bottlenecks, such as the SR-91 connecting the inland empire to beach cities in the Southern California, has motivated many studies on urban transportation analysis. It has long been known that drivers' departure time choice behaviors are a fundamental reason of persistent congestion and, along with route choice behaviors, fundamentally determine the spatial and temporal distribution of a number of morning commuting trips in urban transportation networks. The necessity and importance of including such ``endogenous scheduling" in various the traffic assignment process has been discussed in \citep{Small1992Scheduling}. 

User equilibrium (UE) for drivers' departure time choice at a single bottleneck has been well studied since \citep{vickrey1969congestion,hendrickson1981schedule,small2015bottleneck}.
In such a state, ``the journey costs at all  departure times actually used are equal, and (equal to or) less than those which would be experienced by a single vehicle at any unused time'' \citep{wardrop1952ue}. Here a monetary cost for each trip is designated as a linear function in travel time, time early, time late, and/or toll \citep{Arnott1990bottleneck}
In other words, all commuters have the same total cost, and ``no individual user can improve his total cost by unilaterally changing departure times" \citep{hendrickson1981schedule,mahmassani1984due}. 

Such a departure time user equilibrium can be mathematically formulated as solutions of variational inequalities \citep[e.g.,][]{Friesz1993due,Szeto2004departure}, optimization (linear programming) \citep{iryo2007equivalent}, and linear complementarity problems \citep{akamatsu2015corridor}. However, such formulations are purely phenomenological without offering insights on whether and how the equilibrium can be reached. One promising approach is to model the day-to-day dynamics for departure time choice. However, \citep{iryo2008analysis} demonstrated that Smith's day-to-day dynamical system for departure time choice is unstable and argued that the instability is caused by non-monotonicity of the scheduling cost in time.
\citep{bressan2012variational} proposed two partial differential equations for day-to-day departure flow shifting and  demonstrated the instability of these systems; they also noted that the proof of the instability is challenging with nonlocal operators due to non-monotone scheduling costs, and even the linearized stability analysis was an open problem. 
\citep{guo2016departure,guo2017day} demonstrated that five systems of day-to-day flow shifting dynamics for departure time choices are all unstable and further questioned the existence of a stable user equilibrium at a single bottleneck. \footnote{\citep{benakiva1984dynamic} developed a day-to-day dynamical system model of departure time choice with the logit model, which is different from the UE model. This system is based on a standard flow swapping idea and written in a partial differential equation, but it is not clear whether the system is stable or converges.}

Logically, however, existing studies only examined a finite number of dynamical systems for departure time choice but do not rule out the existence of an unknown stable one. More importantly, empirical observations show that traffic congestion patterns are relatively consistent during peak periods on the same day of a week \citep{jin2015existence}, and our personal driving experience also suggests that the daily departure times for commuting trips are relatively fixed. Therefore, there is an apparent discrepancy between theory and reality regarding drivers' departure time choice.

In this study we attempt to reconcile the discrepancy by presenting a stable day-to-day dynamical system for departure time choice at a single bottleneck. Following \citep{small2015bottleneck} we refer to the departure times from the origin and the arrival times to the destination. They are called the arrival and departure times of the bottleneck, respectively, in many other studies. 
We argue that drivers' departure and arrival time choices are driven by their scheduling payoff choice. 
The biggest difference between existing and this dynamical systems is that, instead of directly modeling the nonlocal departure and arrival time choices, the new model describes the local scheduling payoff choice. 
After presenting three behavioral principles, conceptual models, and mathematical models, we demonstrate that the new system is stable.
Thus the dynamical system model is a behavioral model of user equilibrium.

The rest of the paper is organized as follows. In Section 2, we present definitions related to the within-day dynamics with a point queue model, cost functions, and user equilibrium. In Section 3, we discuss in details the three behavioral rules and corresponding conceptual models. In Section 4, we present mathematical models. In Section 5, we present the corresponding discrete models. In Section 6, we prove that the equilibrium point of the mathematical model leads to the user equilibrium and that it is stable. In Section 7, we present a numerical example to verify the theoretical results. In Section 8, we conclude the study with some potential extensions.

\section{Definitions}

\subsection{Traffic flow variables and point queue model}

We model a bottleneck as a point queue illustrated in \reff{bottleneck_problem}(a), where $F(r,t)$ and $f(r,t)$ are respectively the  departure cumulative flow and  flow-rate at time $t$ on day $r$, $G(r,t)$ and $g(r,t)$ respectively the arrival cumulative flow and flow-rate, $\delta(r,t)$ the queue size,  $\Upsilon(r,t)$ the queueing time for vehicles arriving at $t$, $\Upsilon'(r,t)$ the queueing time for vehicles departing at $t$, and $C$ the bottleneck capacity (maximum service rate) of the queue. The travel demand during a peak period, i.e., the total number of vehicles, is $N$. Here the capacity, $C$, and the travel demand, $N$, are assumed to be constant from day to day. For simplicity we assume that the free-flow travel time is zero, and a non-zero free-flow travel time can be easily incorporated into the discussions.

From the definitions, we have the following relations: 
\bqn
f(r,t)&=&\pd{}t F(r,t),\\
g(r,t)&=&\pd{}t G(r,t),\\
\delta(r,t)&=&F(r,t)-G(r,t),\\
G(r,t)&=&F(r,t-\Upsilon(r,t)), \label{def:fifo}\\
F(r,t)&=&G(r,t+\Upsilon'(r,t)), \label{def:fifo2}
\eqn
where the last two equations are derived from the assumed First-In-First-Out (FIFO) rule.

In addition, from \reff{bottleneck_problem}(a) we can see that 
\bqn
\Upsilon'(r,t)&=&\frac{\delta(r,t)}{C}. \label{def:upsilon_p}
\eqn
Further from $\Upsilon(r,t)=\Upsilon'(r,t-\Upsilon(r,t))$,  the queueing time and the queue size are related by the following equation:
\bqn
\Upsilon(r,t)&=&\frac{\delta(r,t-\Upsilon(r,t))}{C}. \label{def:upsilon}
\eqn
Then we have the following lemma.
\begin{lemma} \label{lemma:zeroqueue}
The following three statements are equivalent at $t$ on day $r$: (i) The queue size is zero; i.e., $\delta(r,t)=0$; (ii) The departure and arrival cumulative flows are equal; i.e., $F(r,t)=G(r,t)$; (iii) The queueing time is zero; i.e., $\Upsilon(r,t)=\Upsilon'(r,t)=0$.
\end{lemma}

\bfg\bc
\includegraphics[width=5.5in]{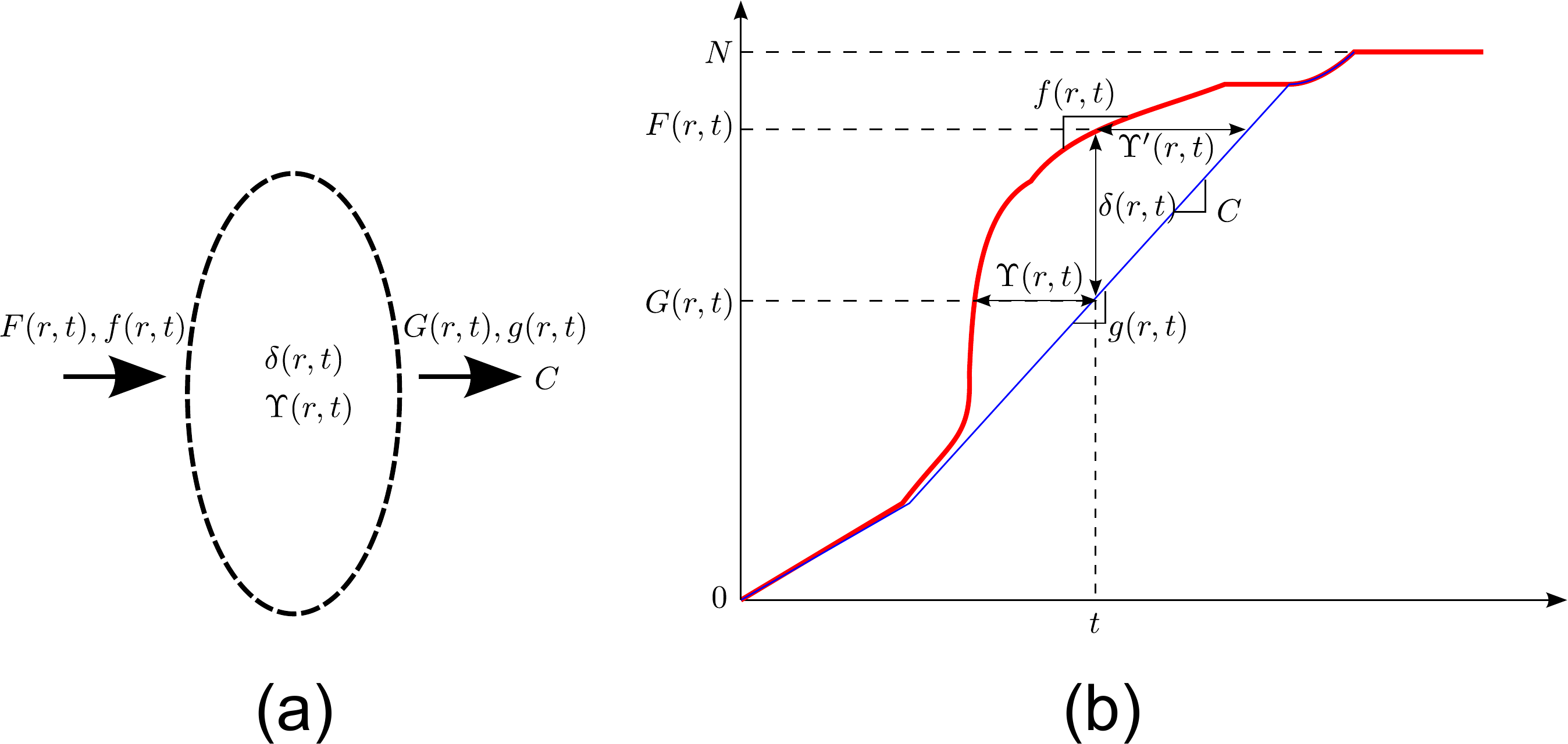}
\caption{A single bottleneck as a point queue and definitions of traffic flow variables on day $r$}
\label{bottleneck_problem}
\ec\efg

Here we apply the point queue model to describe the traffic dynamics \citep{jin2015pqm}
\bqn
\pd{}t \delta(r,t)&=& f(r,t)-g(r,t), \label{pq1}\\
g(r,t)&=&\min\{\frac{\delta(r,t)}{\epsilon}+f(r,t), C \}. \label{arrival-rate}
\eqn
where the storage capacity of the point queue is assumed to be infinite, and $\epsilon=\lim_{\Delta t\to 0^+} \Delta t$ is an infinitesimal hyperreal number and equals $\Delta t$ in the discrete version. 

From \refe{arrival-rate} we have (i) $0\leq g(r,t)\leq C$; i.e., the arrival flow-rate is bounded (but the departure flow-rate may not be); and (ii) $g(r,t)=C$ if $\delta(r,t)>0$. Therefore, if the queue size $\delta(r,t)>0$, the arrival flow-rate has to equal the capacity, and the bottleneck is fully utilized.
Equivalently, if the arrival flow-rate is smaller than the capacity; i.e., if the bottleneck is under-utilized, then the queue size $\delta(r,t)=0$. Note that, when the arrival flow-rate equals the capacity, it is possible that the queue size is zero; and vice versa.  Thus we have
\bqn
\delta(r,t) (g(r,t)-C)&=&0, \label{queue-arrival-complementariy-1}
\eqn
which is a complementarity condition between the queue  size and the arrival flow-rate and can be re-written as  $\delta(r,t)\geq 0$, $g(r,t) \leq C$, and $\min\{\delta(r,t), C-g(r,t) \}=0$.

From Lemma \ref{lemma:zeroqueue}, \refe{queue-arrival-complementariy-1} is equivalent to the following conditions \citep{iryo2007equivalent}:   
\bqn
(F(r,t)-G(r,t)) (g(r,t)-C)&=&0, \\
\Upsilon(r,t) (g(r,t)-C)&=&0,  \label{queue-arrival-complementariy}\\
\Upsilon'(r,t) (g(r,t)-C)&=&0.
\eqn

From the above discussions we have the following lemma.

\begin{lemma}\label{lemma:inverse}
	If $g(r,t)<C$; i.e., if the bottleneck is under-utilized at $t$, then $\delta(r,t)=0$,  $F(r,t)=G(r,t)$, $f(r,t)=g(r,t)$, and $\Upsilon(r,t)=\Upsilon'(r,t)=0$. Equivalently, when $\delta(r,t)>0$, $F(r,t)>G(r,t)$, or $\Upsilon(r,t)>0$, then $g(r,t)=C$.
\end{lemma}

\refe{pq1} can be re-written as
\bqn
\pd{}t \delta(r,t)&=&\max\{-\frac{\delta(r,t)}{\epsilon}, f(r,t)-C \}, \label{queue_rate}
\eqn
whose discrete version is
\bqn
\delta(r,t+\dt)&=&\max\{0, \delta(r,t)+(f(r,t)-C) \dt\}. \label{discrete-queue}
\eqn

\begin{lemma}\label{departure_arrival_lemma}
	If $f(r,t)>0$, then $g(r, t+\Upsilon'(r,t))>0$; i.e., if there are vehicles departing at $t$, then there are vehicles arriving at $t+\Upsilon'(r,t)$. 
	Similarly, if $g(r,t)>0$, then $f(r,t-\Upsilon(r,t))>0$; i.e., there are vehicles arriving at $t$, then there are vehicles departing at $t-\Upsilon(r,t)$. Note that if $f(r,t)=0$, $g(r, t+\Upsilon'(r,t))$ could be either 0 or positive; and if $g(r,t)=0$, $f(r,t-\Upsilon(r,t))$ could also be either 0 or positive.
	\end{lemma}
{\em Proof}. From \refe{def:fifo2}, \refe{def:upsilon_p},  and \refe{queue_rate}, we have
\bqs
f(r,t)&=&g(r,t+\Upsilon'(r,t)) \cdot (1+\frac{\pd{}t \delta(r,t)}C)=g(r,t+\Upsilon'(r,t)) \cdot \max\{ 1-\frac{\delta(r,t)}{C\epsilon}, \frac{f(r,t)}C \}.
\eqs
Thus if $f(r,t)>0$, then $\max\{ 1-\frac{\delta(r,t)}{C\epsilon}, \frac{f(r,t)}C \}>0$, and $g(r,t+\Upsilon'(r,t))>0$.

From \refe{def:fifo}, \refe{def:upsilon}, and \refe{queue_rate}, we have
\bqs
g(r,t)&=&f(r,t-\Upsilon(r,t)) \cdot (1-\pd{}t \Upsilon(r,t))=f(r,t-\Upsilon(r,t)) \cdot \frac{C}{C+\pd{}t \delta(r,t-\Upsilon(r,t))}\\
&=& f(r,t-\Upsilon(r,t)) \cdot \frac{C}{\max\{C-\frac{\delta(r,t-\Upsilon(r,t))}{\epsilon}, f(r,t-\Upsilon(r,t)) \}}
.
\eqs
Thus if $g(r,t)>0$, then $f(r,t-\Upsilon(r,t))>0$, assuming that the queueing time is properly defined when the departure flow-rate is zero.
\eop

\subsection{Costs and user equilibrium}
A driver's cost at a single bottleneck comprises of two parts: the queueing cost, caused by congestion, and  the scheduling cost, caused by the schedule delay.

The queueing cost for a driver arriving at $t$ is \citep{arnott1990departure}
\bqn
\phi_1(r,t)&=&\alpha \Upsilon(r,t), \label{def:cost1}
\eqn
where $\alpha=$\$6.4/hr. Here we assume that all drivers have the same free-flow travel time through the point queue, which is zero in this case.

For the scheduling cost, we assume that all vehicles have the same ideal arrival time, $t_*$. A piece-wise linear scheduling cost function for a driver arriving at $t$ can be written as:
\bqn
\phi_2(t)&=&\beta \max\{ t_*-t,0\} +\gamma \max\{t-t_*, 0\}, \label{def:cost2}
\eqn
where the first term on the right-hand side is the penalty for an early arrival, the second term for a late arrival. In \citep{arnott1990departure}, $\beta=$\$3.90/hr, and $\gamma=$15.21/hr.
A necessary condition for the existence of user equilibrium is $\beta<\alpha$; in general, $\beta/\alpha=0.5$, and $\gamma/\alpha=2$ \citep{small2015bottleneck}.

Thus the total cost is
\bqn
\phi(r,t)&=&\phi_1(r,t)+\phi_2(t)=\alpha \Upsilon(r,t)+\beta \max\{ t_*-t,0\} +\gamma \max\{t-t_*, 0\}. \label{def:totalcost}
\eqn
The definitions of costs are illustrated in \reff{bottleneck_costs}, where the height of the shaded areas represent the arrival flow-rates. Here we assume that the coefficients, $\alpha$, $\beta$, and $\gamma$, are constant from day to day. Thus the scheduling cost $\phi_2(t)$ is independent of $r$, but the queueing cost $\phi_1(r,t)$ is not. In addition, we assume that all drivers have the same coefficients; thus all drivers belong to the same class.

\bfg\bc
\includegraphics[width=5.5in]{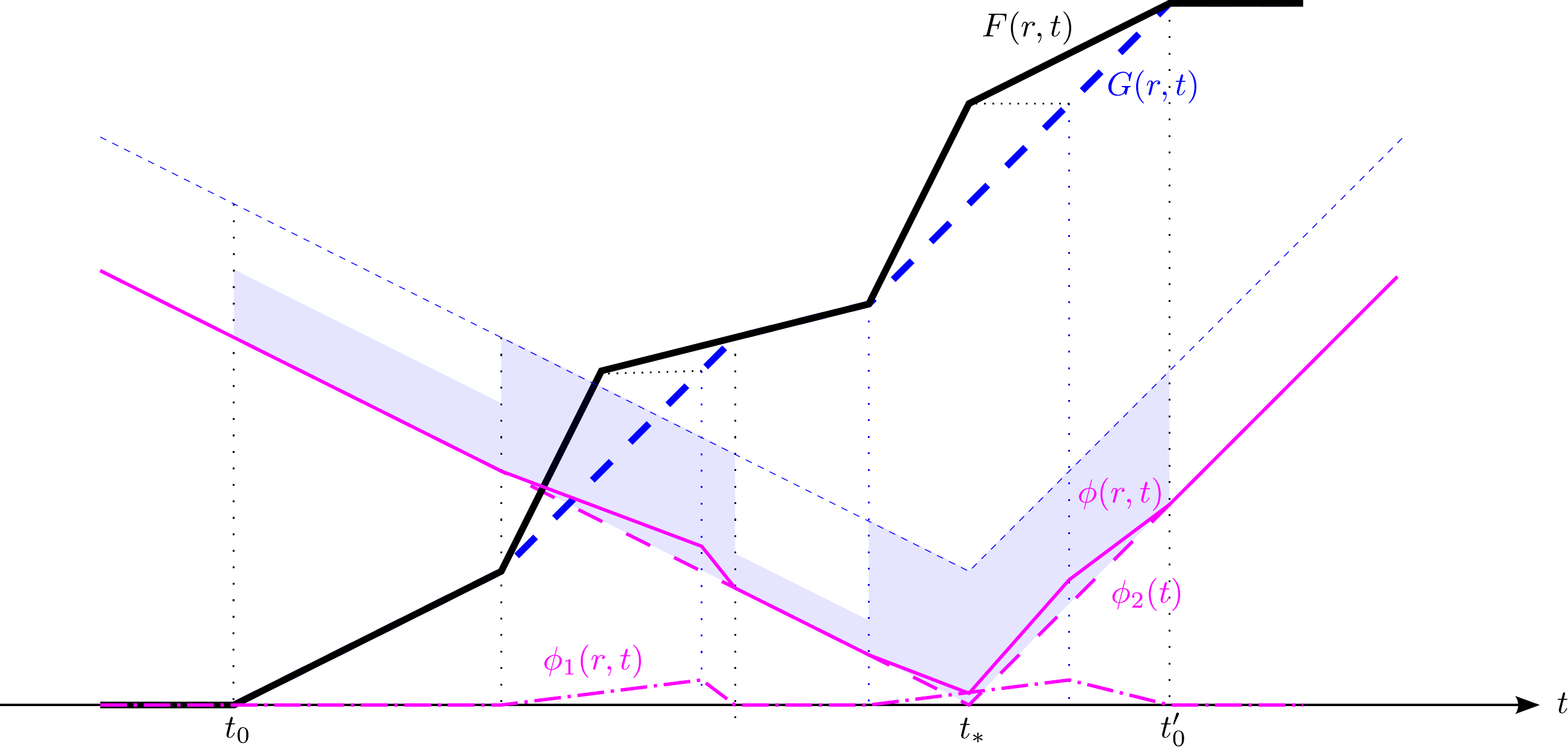}
\caption{Illustration of travel, scheduling, and total costs: $\alpha=\frac 12$, $\beta=\frac 14$, and $\gamma=1$}
\label{bottleneck_costs}
\ec\efg

\begin{lemma} \label{lemma:vacancy}
If the arrival flow-rate at a time is smaller than the capacity; i.e., if the bottleneck is under-utilized, then the total cost equals the scheduling  cost.
\end{lemma}
{\em Proof}. From \refe{queue-arrival-complementariy}, if the arrival flow-rate is smaller than the capacity, then the queueing time has to be zero. Then from \refe{def:totalcost}, the total cost equals the scheduling cost. \eop

\begin{theorem} \label{theorem:vacancy}
	If the bottleneck is under-utilized at $t_1$ and the scheduling cost at $t_2$ is not smaller than that at $t_1$, then the total cost at $t_2$ is not smaller than that at $t_1$. That is, if $g(r,t_1)<C$, and $\phi_2(t_2)\geq\phi_2(t_1)$, then $\phi(r,t_2) \geq \phi(r,t_1)$.
\end{theorem}
{\em Proof}. From Lemma \ref{lemma:vacancy}, we have $\phi(r,t_1)=\phi_2(t_1)$. Thus
\bqs
\phi(r,t_2)&\geq&\phi_2(t_2) \geq \phi_2(t_1)=\phi(r,t_1).
\eqs
That is, the total cost at $t_2$ is not smaller either. 
Note that this conclusion may not be true if the bottleneck at $t_1$ is fully utilized.
\eop

\begin{definition} \label{def:DTUE:ATUE}
	A system reaches the {\bf departure time user equilibrium (DTUE)} if all used departure times have the same total cost, which is not greater than those of unused departure times. That is, $\phi(r,t+\Upsilon'(r,t))=\phi^*$ if $f(r,t)>0$; and $\phi(r,t+\Upsilon'(r,t))\geq \phi^*$ if $f(r,t)=0$. Here $\phi^*$ is the minimum cost.
	
In contrast, a system reaches the {\bf arrival time user equilibrium (ATUE)} if all used arrival times have the same total cost, which is not greater than those of unused arrival times. That is, $\phi(r,t)=\phi^*$ if $g(r,t)>0$; and $\phi(r,t)\geq \phi^*$ if $g(r,t)=0$.		
	\end{definition}

\begin{theorem} \label{theorem:DTUE-ATUE-equivalence}
	The DTUE and ATUE are equivalent.
	\end{theorem}
{\em Proof}. Assume that a system is in the ATUE state.  (i) If $f(r,t)>0$, then $g(r,t+\Upsilon'(r,t))>0$ from Lemma \ref{departure_arrival_lemma}. From the ATUE definition, we have $\phi(r,t+\Upsilon'(r,t))=\phi^*$. (ii) If $f(r,t)=0$, then $g(r,t+\Upsilon'(r,t))\geq 0$. From the ATUE definition, we have $\phi(r,t+\Upsilon'(r,t))\geq\phi^*$. Thus the system is the DTUE state.

Similarly we can use Lemma \ref{departure_arrival_lemma} to prove that DTUE also implies ATUE. Therefore the two equilibrium states are equivalent.  \eop

\section{Behavioral principles and conceptual models for departure time choice}
If initially all drivers attempt to minimize their scheduling costs and depart at the same ideal time without considering congestion, the departure flow-rate is effectively infinite, the bottleneck becomes congested, and all except the first driver will experience additional congestion costs. Thus drivers' departure time choice behaviors have to be much more sophisticated in order to minimize their actual costs. 

 Given the complete information of the current day, including the departure and arrival cumulative flows as well as the travel, scheduling, and total costs as shown in \reff{bottleneck_costs},\footnote{In reality drivers may accumulate such information by learning from various traffic information sources, others' experience, or their own trials.}  drivers can choose both their departure times and arrival times. Even though drivers depart from the origin before arriving at the destination, in our daily decision process ``an individual with a fixed destination appointment (usually) works backward from that appointment time to determine which scheduled public transportation vehicle or train he must use in order to arrive no later than his appointment time''  \citep{kraft1967new}. That is, drivers are mostly concerned about their schedule delays, which equal the differences between the actual and ideal arrival times. Hence the first behavioral principle regarding departure time choice is that drivers actually choose their departure and arrival times in a backward fashion; i.e., they determine arrival times first and then choice the departure times accordingly. We refer to this as the backward choice principle.

\subsection{Departure time choice}

It has been observed in \citep{vickrey1969congestion} that ``those arriving at their offices closest to their desired times generally have to spend relatively more time in the queue than those who choose to push their arrival time further away from the desired time'', and in \citep{hendrickson1981schedule} that ``for later arrivals at the bottleneck, the queue must increase so that the reduction in schedule delay is exactly balanced by an increase in the queuing time.'' 
Thus the second behavioral principle is that drivers update their departure times to balance the total costs; i.e., their queueing costs have to be larger if their scheduling costs are smaller. We refer to this as the cost balancing principle.

\commentout{
The cost balancing principle obviously does not apply to those who arrive during an underutilized period, since the corresponding travel (queueing) cost is zero. Therefore, the time axis contains multiple fully-utilized intervals, during which the arrival flow-rates equal the capacity. 

Given an arrival cumulative curve, a feasible departure cumulative flow should satisfy the following conditions: (i) The slope of $F(t)$ can be infinite (when all vehicles depart at the same time), but not negative (which leads to a negative departure flow-rate); (ii) The departure cumulative flow cannot be smaller than the arrival cumulative flow; and (iii) The queueing cost at the beginning of the interval has to be zero, but that at the end of the interval may not be zero, when the last vehicle departs at the same time with other vehicles.  Clearly there can be many feasible departure cumulative flows: the minimum one is the same as the arrival cumulative flow when all drivers choose the same departure times as the arrival times; but in this case the total cost is least balanced.}

\bfg\bc
\includegraphics[width=5.5in]{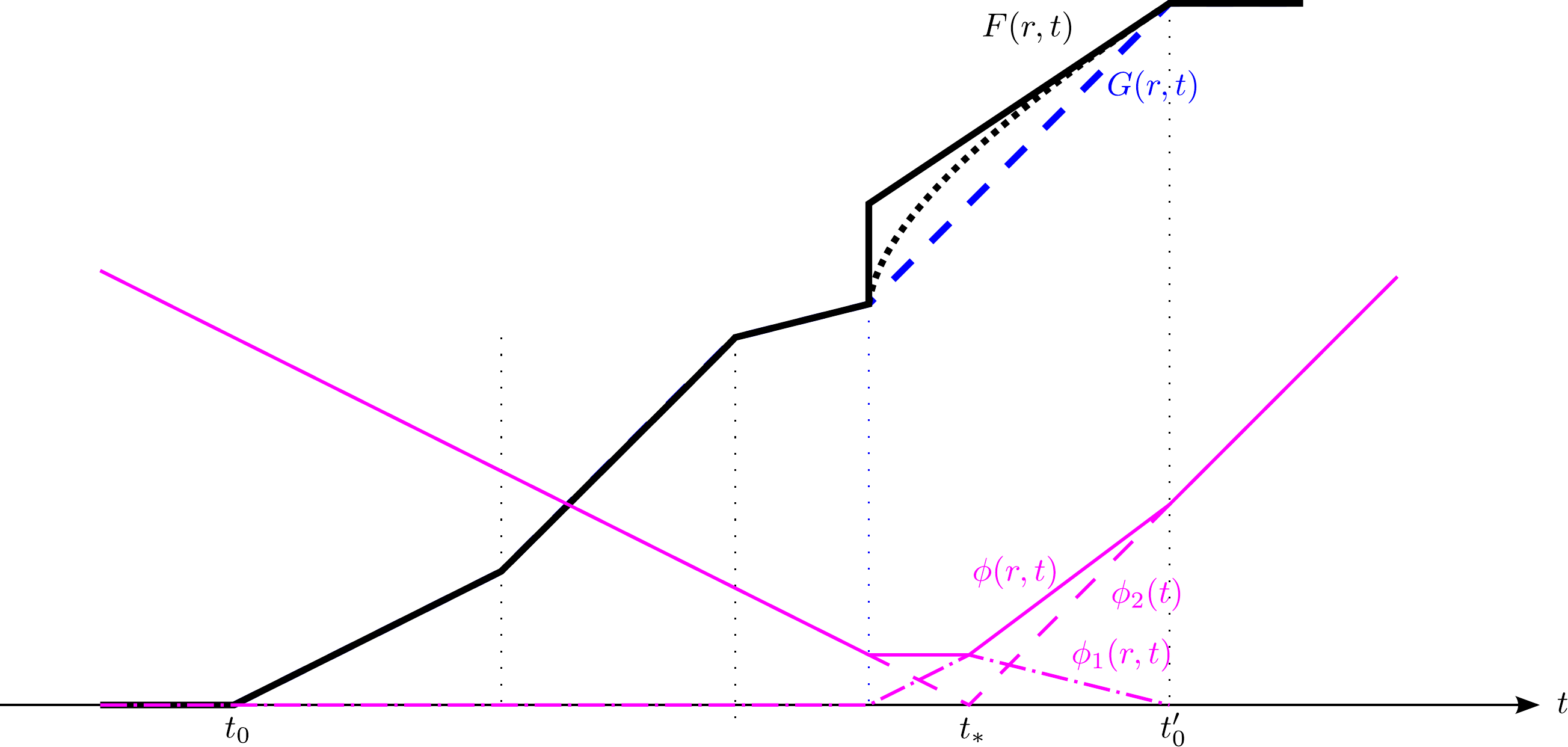}
\caption{The cost balancing principle}
\label{bottleneck_balance}
\ec\efg

Note that no queue exists to balance an under-utilized time, since the departure flow-rate automatically equals the arrival flow-rate. We can choose a departure cumulative flow to balance the total cost during each fully-utilized interval either locally within the interval or globally across multiple intervals. Here we propose the following simple conceptual model, as illustrated by \reff{bottleneck_balance}, where  the arrival cumulative flow  is given by the thick dashed curve, the partially balanced departure cumulative flow  by the thick solid curve, and the corresponding total, scheduling, and traveling costs by the solid, dashed, and dash-dotted curves at the bottom.  (i) For an interval before the ideal arrival time, all drivers choose the same departure times as the arrival times, such that their total costs are minimum without any queueing cost. (ii) For an interval after the ideal arrival time, all drivers choose the same departure times as the arrival times, such that their total costs are minimum without any queueing cost. (iii) For an interval containing the ideal arrival time, vehicles arriving before the ideal arrival time choose their departure times such that their total costs equal the maximum total cost (the first vehicle's cost), and vehicles arriving later choose their departure times such that their total costs form a continuous line; in particular, if the two boundaries have the same scheduling cost, then all vehicles leaving in this interval have the same cost, and the solid curve at the bottom of the figure is flat in this interval. In this model, the drivers' departure time choice is totally determined by their arrival time choice, and the departure time choice dynamics predicate on the arrival time choice dynamics, confirming the first behavioral principle.
In \reff{bottleneck_balance}, other curves are also feasible (one example is the thick dotted curve), but here we choose the simple cost balancing principle.

\subsection{Arrival time choice}
The third behavioral principle is that a driver chooses his/her arrival time to reduce the scheduling cost. In particular, he tries to switch his arrival time to  an under-utilized arrival time with a smaller scheduling cost.
According to Theorem \ref{theorem:vacancy}, such a switch reduces not only  the scheduling cost, but also the total cost.  We refer to this as the scheduling cost reducing principle. Note that switching to a fully-utilized arrival time with a smaller scheduling cost may not be beneficial.

\bfg\bc
\includegraphics[width=5.5in]{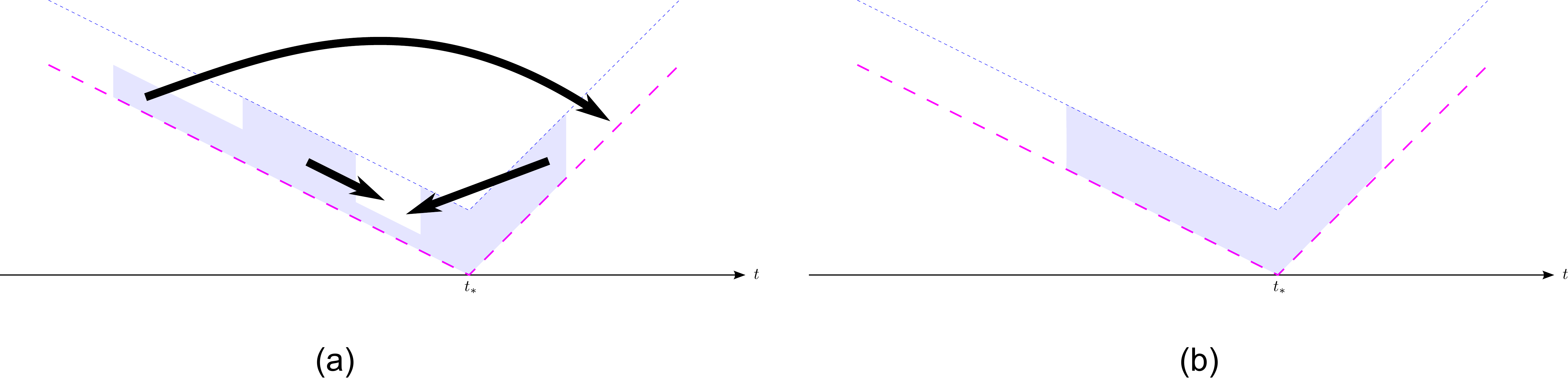}\caption{A V-shaped tube model of the scheduling cost reducing principle for arrival time choice: (a) Initial condition; (b) Equilibrium condition}
\label{bottleneck_tube}
\ec\efg

In \reff{bottleneck_tube},  the heights of the dashed curves represent the scheduling costs at different times, and the vertical gaps between the dashed and dotted curves represent the bottleneck capacity. \reff{bottleneck_tube}(a) shows an example of the initial arrival flow-rates: the height of the shaded region at a time represents the arrival flow-rate. If we view the region bounded by the dashed and dotted curves as a V-shaped tube, then the figure clearly shows when the bottleneck is under-utilized or not: the bottleneck is fully utilized when the tube is filled by the shaded region, and underutilized when there is vacancy in the tube.

As shown in \reff{bottleneck_tube}(a), the three arrows represent three potential switches of arrival times that satisfy the scheduling cost reducing principle, where the driver with the arrival time at an arrow tail switches his/her arrival time to that at the arrow head. When the system reaches the ATUE, nobody can find a vacancy inside the V-shaped tube with a smaller scheduling cost, as shown in \reff{bottleneck_tube}(b). In this sense, the third principle of arrival time choice is equivalent to finding a lower vacancy  in the V-shaped tube. Note that, as shown in the figure, due to the non-monotonicity in the scheduling cost function, drivers can shift to earlier or later arrival times to reduce their scheduling costs. That is, such arrival time shifting can be nonlocal. Such a nonlocal process also occurs for departure time shifting \citep{bressan2012variational}.

 Conceptually such movements are consistent with the movements of fluids (e.g. water) caused by gravity in a V-shaped tube: they will attempt to fill a vacancy at a lower point, and the heights of both sides are the same in  equilibrium.  
However, we cannot directly use fluid dynamics equations to describe the day-to-day arrival time shifting dynamics, since (i) the equilibrium surface of the fluid is not flat in \reff{bottleneck_tube}(b); (ii) the former can lead to oscillations at the equilibrium in a frictionless V-shaped tube \citep{landau1987fluid}, and (iii) more importantly, the latter involves nonlocal behaviors, which cannot be captured in the former. 

\subsection{Scheduling payoff choice}
To enable the nonlocal movement of fluids in a V-shaped tube, we can connect the two sides of the V-shaped tube with capillaries with negligible widths, as shown in \reff{single_tube}(a). Then the V-shaped tube can be simplified into a single tube as shown in \reff{single_tube}(b). 

We introduce a new variable $x$ for the scheduling payoff, equal to the negative scheduling cost, as shown in \reff{single_tube}(b): 
\bqn
x&=&-\phi_2(t).
\eqn
Therefore, a value of $x$ corresponds to two arrival times, $t_1(x)$ and $t_2(x)$, as shown in \reff{single_tube}(a), where
\bsq\label{def:t1t2}
\bqn
x&=&-\beta (t_*-t_1(x)),\\
x&=&-\gamma(t_2(x)-t_*).
\eqn
\esq

Therefore, the third behavioral principle for arrival time choice can be directly applied to scheduling payoff choice: all drives attempt to gain their scheduling payoff from day to day. Thus the scheduling cost reducing principle can also be called the scheduling payoff gaining principle. Again, such movements are consistent with those of the fluid caused by gravity in a single tube: the fluid in the single tube moves to the bottom part; i.e., the fluid  moves in the positive direction of $x$.
Thus the nonlocal arrival time choice problem in a V-shaped tube is equivalent to the local scheduling payoff choice problem in a single tube. 
In a sense, the single tube can be considered an imaginary unidirectional road, where vehicles move toward larger scheduling payoff or smaller scheduling cost.

\bfg\bc
\includegraphics[width=5.5in]{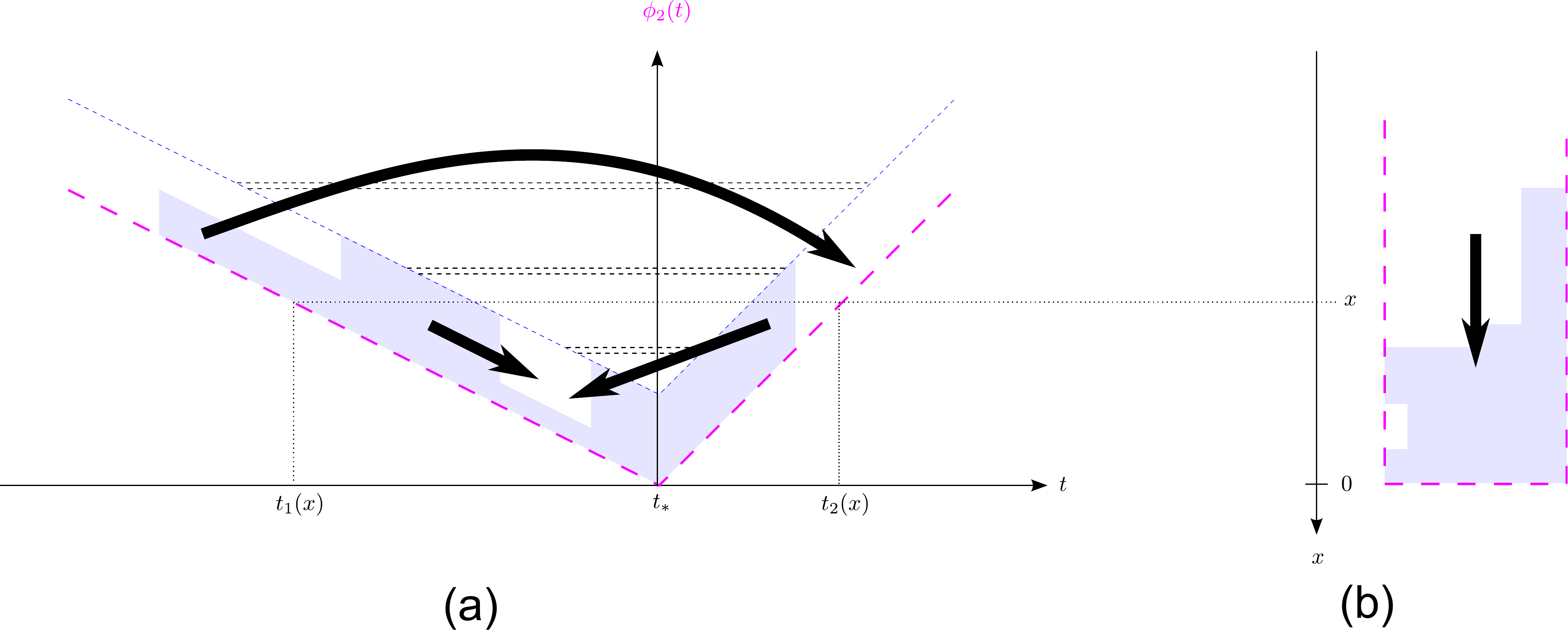}\caption{A single-tube model of scheduling payoff choice} \label{single_tube}
\ec\efg

\section{Mathematical models}
In this section we mathematically formulate the conceptual models in the preceding section. We present the continuous models of scheduling payoff choice and  corresponding arrival and departure flow-rates.

We introduce a new variable, $k(r,x)$, whose unit is veh/\$, for the number of vehicles per dollar choosing the scheduling payoff of $x$ on day $r$. It is referred to as the imaginary density as it measures the number of vehicles per unit length on the imaginary road. Thus the total number of vehicles with scheduling payoffs not greater than $x$ is given by 
\bqn
\int_{y=x}^0 k(r,y) dy&=& \int_{t_1(x)}^{t_2(x)} g(r,t) dt, \label{def:rho_int}
\eqn
which is the integral form of the relation between $k(r,x)$ and $g(r,t)$. The corresponding differential form is 
\bqn
k(r,x)&=& \frac 1 \beta g(r,t_1(x))  + \frac 1 \gamma g(r,t_2(x)). \label{def:rho}
\eqn
That is, the vehicle density at $x$ equals the sum of the arrival flow-rates at the two time, $t_1(x)$ and $t_2(x)$, adjusted by $\beta$ and $\gamma$ respectively. It is represented by the horizontal width at $x$ of the shaded region in \reff{single_tube}(b).

Correspondingly, the width of the single tube at $x$ shown in \reff{single_tube}(b) equals the sum of the capacities at the two arrival times, adjusted by $\beta$ and $\gamma$ respectively. We refer to the width as the ``imaginary jam density", which is denoted by $\kappa$; thus
\bqn
\kappa=(\frac 1\beta+\frac 1\gamma)C. \label{def:kappa}
\eqn

\subsection{The LWR model of scheduling payoff choice}

The scheduling payoff choice problem shown in \reff{single_tube}(b) can be viewed as a fluid flow problem, in which $k(r,x)$ is the fluid density. Furthermore, from
\bqn
\int_{x=-\infty}^0 k(r,x) dx&=&\int_{t=-\infty}^{\infty} g(r,t) dt=N,
\eqn
we can see that the imaginary density is conserved.
Also as shown in \reff{single_tube}(b), the fluid always travels toward $x=0$. 
For such a unidirectional, conservative fluid flow system we can apply the Lighthill-Whitham-Richards (LWR) model to describe its dynamics \citep{lighthill1955lwr,richards1956lwr}:
\bqn
\pd{}r k(r,x)+ \pd{}x Q(k(r,x))&=&0, \label{lwr_sc}
\eqn
subject to 
\ben
\item $k(r,x) \in[0,\kappa]$; 
\item The fundamental diagram, e.g., the triangular fundamental diagram 
\bqn
Q(k)&=&\min\{u k, w(\kappa-k) \}, \label{lwr_tri_fd}
\eqn
where $u$ is the free-flow speed and $w$ the shock wave speed in congested traffic. Here both $u$ and $w$ are positive and can choose arbitrary values. The unit of $Q(k)$ is veh/day, and that of $u$ and $w$ is \$/day.
\een

Following the kinematic wave theory based on the Cell Transmission Model \citep{daganzo1995ctm,lebacque1996godunov,jin2009sd}, we introduce the demand and supply variables as $d(r,x)$ and $s(r,x)$, where 
\bsq
\bqn
d(r,x)&=&D(k(r,x))\equiv u \min\{k(r,x), \kappa_c\},\\
s(r,x)&=&S(k(r,x))\equiv w (\kappa -\max\{k(r,x),\kappa_c \} ),
\eqn
where $\kappa_c=\frac{w}{u+w} \kappa$ is the critical density.
\esq
In this sense, the boundary condition of \refe{lwr_sc} is given by 
\bsq\label{bc_sp}
\bqn
s(r, 0^+)&=&0; 
\eqn
i.e., the downstream boundary is closed, and the fluid cannot move beyond the bottom point of the single tube. 
Without loss of generality, we assume that $\phi_2(t_0)=\phi_2(t_0')=L$, which is the maximum scheduling cost. Correspondingly, the range of the scheduling payoff is $x\in [-L,0]$. Since no drivers choose a departure time beyond $[t_0,t_0']$, $k(r,x)=0$ for $x\notin [-L,0]$. As the fluid flow travels in the direction of $x$, the external demand is also zero:
\bqn
d(r,-L^-)=0.
\eqn
\esq

On day $0$, given the initial departure cumulative flow, $F(0,t)$, we can use the point queue model in Section 2 to calculate $G(0,t)$, $g(0,t)$, and $k(0,x)$, which is the initial condition for \refe{lwr_sc}. 

The LWR model, \refe{lwr_sc}, subject to the initial density $k(0,x)$ and the boundary conditions, \refe{bc_sp}, is an infinite-dimensional dynamical system model, in which the imaginary density at different payoffs evolve with respect to the day variable, $r$. This is a typical link initial-boundary value problem (LIBVP) of the LWR model, and its solution describes the day-to-day dynamics of scheduling payoff choice on the imaginary road. 

\subsection{Evolution of the arrival and departure flow-rates} \label{section:costbalancing}

With the day-to-day dynamics of scheduling payoff choice described by the LWR model, \refe{lwr_sc}, we first update the day-to-day evolution of the arrival flow-rate.
From \refe{def:rho}, we can have different arrival flow-rates at $t_1(x)$ and $t_2(x)$ with a given $k(r,x)$. Here we assume an equal splitting principle:  
\bqn
g(r,t_1(x))&=&g(r,t_2(x))=\frac {\beta\gamma} {\beta+\gamma} k(r,x), \label{density2arrival}
\eqn
which are between $0$ and $C$.

\begin{lemma} \label{lemma:xtog}
	$k(r,x)=\kappa$ if and only if $g(r,t_1(x))=g(r,t_2(x))=C$. Similarly,	$k(r,x)=0$ if and only if $g(r,t_1(x))=g(r,t_2(x))=0$. Furthermore, $k(r,x)\in[0,\kappa]$ if and only if $g(r,t)\in[0,C]$.
	\end{lemma}
{\em Proof}. Note that $t_1(x)$ and $t_2(x)$ are defined in \refe{def:t1t2}. From \refe{density2arrival} and that $k(r,x)\in [0, \kappa]$ with $\kappa$ given in \refe{def:kappa}, we can simply prove the lemma. \eop

On day $r$, we denote $(x_*^r,0]$ as the largest jammed scheduling payoff interval containing $0$. That is, $k(r,x)=\kappa$ for $x\in (x_*^r,0]$; but $k(r,x)<\kappa$ for a slightly smaller $x=x_*^r-\epsilon$, where $\epsilon$ is an infinitesimal positive number. Note that the interval is empty if $x_*^r=0$. Thus we can divide the study time period $[t_0,t_0']$ into three intervals: $[t_0, t_1(x_*^r)]$, $(t_1(x_*^r), t_2(x_*^r))$, and $[t_2(x_*^r), t_0']$, where the two intermediate times $t_1(x_*^r)$ and $t_2(x_*^r)$ satisfy the following conditions: (i) $\phi_2(t_1(x_*^r))=\phi_2(t_2(x_*^r))=-x_*^r$; (ii) $t_0<t_1(x_*^r)\leq t_* \leq t_2(x_*^r)<t_0'$; (iii) $g(r,t)=C$ for $t\in (t_1(x_*^r), t_2(x_*^r))$; (iv) $g(r,t)<C$ for $t=t_1(x_*^r)-\epsilon$ or $t=t_2(x_*^r)+\epsilon$ with $\epsilon$ as an infinitesimal positive number. That is, the second interval contains $t_*$, but may be of zero length when  $t_1(x_*^r)= t_* = t_2(x_*^r)$; it is fully utilized, but the times immediately outside the interval are under-utilized. In other words, the second interval is the largest fully utilized interval containing the ideal arrival time.\footnote{Without loss of generality, we assume that $t_0$ and $t_0'$ are far from the ideal arrival time such that the first and last intervals are always non-empty.}
 
 For under-utilized sub-intervals in the first and last time intervals, $f(r,t)=g(r,t)<C$ and $F(r,t)=G(r,t)$, according to Lemma \ref{lemma:inverse}; for fully utilized sub-intervals, the departure time costs can never be balanced inside these sub-intervals, and we simply set $f(r,t)=g(r,t)=C$ and $F(r,t)=G(r,t)$.
 If the length of the second interval $(t_1(x_*^r), t_2(x_*^r))$ is positive; i.e., if $t_1(x_*^r)<t_*<t_2(x_*^r)$, 
 we set with the cost balancing principle $\phi(r,t)=-x_*^r$ for $t\in (t_1(x_*^r), t_2(x_*^r))$. There are the following two sub-cases: (i) for $t\in(t_1(x_*^r),t_*]$, we have
 \bqs 
 \alpha \Upsilon(r,t)+\beta (t_*-t)&=&\beta (t_*- t_1(x_*^r))=-x_*^r;
 \eqs
for $t\in (t_*, t_2(x_*^r))$, we have
\bqs
\alpha \Upsilon(r,t)+\gamma(t-t_*)&=&\gamma (t_2(x_*^r)-t_*)=-x_*^r. 
\eqs 
Therefore, the queue length is
 \bqn
 \Upsilon(r,t)&=&\cas{{ll}\frac \beta \alpha (t-t_1(x_*^r)), & t\in(t_1(x_*^r),t_*]; \\
 \frac \gamma \alpha (t_2(x_*^r)-t),	& t\in (t_*, t_2(x_*^r)).}
 \eqn

In the second time interval, the arrival cumulative flow is linear: $G(r,t)=G(r,t_1(x_*^r))+  C (t-t_1(x_*^r)) = G(r,t_2(x_*^r)) +  C (t-t_2(x_*^r))$; in addition, $F(r,t_1(x_*^r))=G(r,t_1(x_*^r))$, and $F(r, t_2(x_*^r))=G(r,t_2(x_*^r))$. Thus we have from \refe{def:fifo} that 
\bqs
F(r,t-\frac \beta \alpha (t-t_1(x_*^r)))&=&F(r,t_1(x_*^r))+C(t-t_1(x_*^r)), \quad t\in(t_1(x_*^r),t_*]; \\
F(r,t-\frac \gamma \alpha (t_2(x_*^r)-t))&=&F(r,t_2(x_*^r))+C(t-t_2(x_*^r)), \quad t\in (t_*, t_2(x_*^r)),
\eqs
which lead to the following cumulative departure flow
\bqn
F(r,t)&=&\min\{F(r,t_1(x_*^r))+\frac {C}{1-\frac \beta \alpha} (t-t_1(x_*^r)), F(r,t_2(x_*^r))+\frac {C}{1+\frac \gamma \alpha} (t-t_2(x_*^r))\},
\eqn
and the following departure flow-rate
\bqn
f(r,t)&=&\cas{{ll} \frac {C}{1-\frac \beta \alpha}, & t\in [t_1(x_*^r), \frac \beta\alpha t_1(x_*^r)+(1-\frac \beta\alpha) t_*);\\
	\frac {C}{1+\frac \gamma \alpha}, & t\in [- \frac \gamma\alpha t_2(x_*^r)+(1+\frac \gamma\alpha) t_*, t_2(x_*^r)).
}	
\label{equilibrated_departure_flowrate}
\eqn
Note that the departure flow-rate and flow are well-defined if and only if $\alpha>\beta$ \citep{small2015bottleneck}.
In addition,
\bqs
\frac \beta\alpha t_1(x_*^r)+(1-\frac \beta\alpha) t_*&=&-\frac \gamma\alpha t_2(x_*^r)+(1+\frac \gamma\alpha) t_*,
\eqs
since $\phi_2(t_1(x_*^r))=\phi_2(t_2(x_*^r))$, or equivalently $\beta (t_*-t_1(x_*^r))=\gamma (t_2(x_*^r)-t_*)$. We denote this time by $t_*^r$. Thus we have at $t_*^r$
\bqs
F(r,t_*^r)&=&F(r,t_1(x_*^r))+C(t_*-t_1(x_*^r))=F(r,t_2(x_*^r))+C(t_*-t_2(x_*^r))=G(r, t_*);
\eqs
that is, $t_*^r$ is the departure time for vehicles with the ideal arrival time. Note that this departure time depends on $x_*^r$.

Even though $\phi(r,t)=-x_*^r$ for $t\in(t_1(x_*^r),t_2(x_*^r)$, and $-x_*^r \leq \phi(r,t)$ for $t\notin(t_1(x_*^r),t_2(x_*^r)$. The departure and arrival flow-rates are only partially equilibrated, if $g(r,t)>0$ for $t\notin(t_1(x_*^r),t_2(x_*^r)$.

\section{Discrete models}\label{section:discrete_models}
\bfg\bc
\includegraphics[width=5.5in]{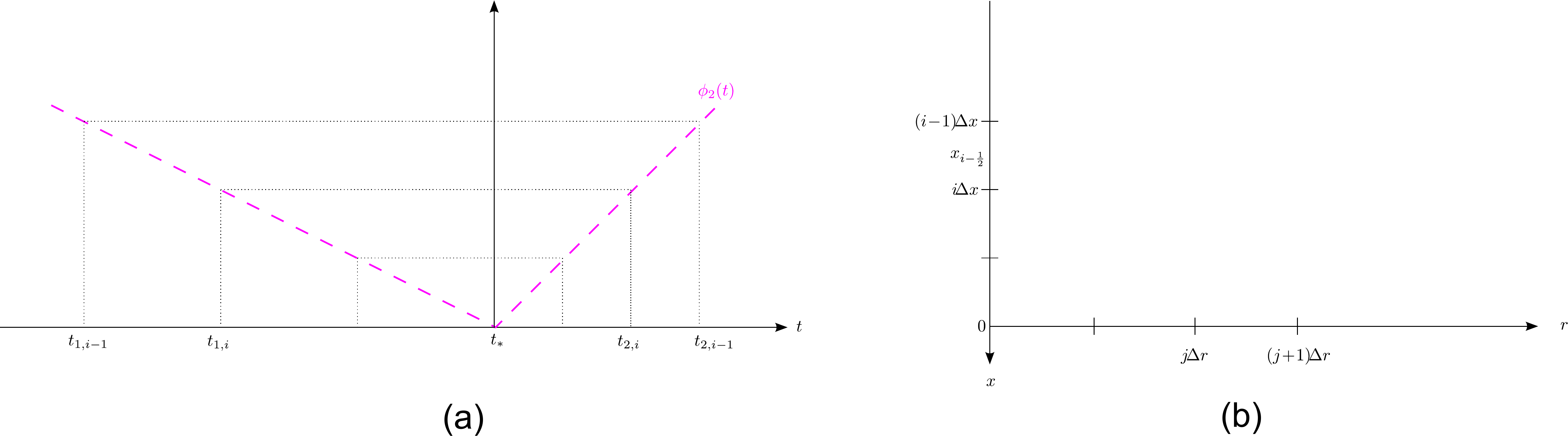}\caption{Discretization of $r$, $x$, and $t$-axes}\label{CTM_sc}
\ec\efg

As shown in \reff{CTM_sc}(b), we discretize the $r$-axis into a number of steps with a day-step size $\Delta r$ and partition the imaginary $x$-axis into $I$ cells with a cell size of $\dx$. The cell boundaries are at $i\dx$ with $i=0,-1,-2,\cdots,-I$. The center of cell $i$ is at $x_{i-\frac 12}=(i-\frac 1 2)\dx$ ($i=0,\cdots, -I+1$). Correspondingly, the $t$-axis is divided into a number of intervals, with the boundaries $t_{1,i}=t_1(i\dx)=t_*+i \frac{\dx}{\beta}$ and $t_{2,i}=t_2(i\dx)=t_*-i \frac{\dx}{\gamma}$. We divide the study period $[t_0,t_0']$ into $M$ intervals; thus $\dt=\frac {t_0'-t_0}{M}$. Here we assume that $\dt$ is generally much smaller than $\frac{\dx}{\beta}$ and $\frac{\dx}{\gamma}$. We denote the center of time interval $m$ by $t_{m-\frac 12}=t_0+(m-\frac 12) \dt$ ($m=1,\cdots,M$).  Then the dynamical system of departure time choice can be numerically solved in the following steps.

\subsection{The initial condition}\label{section:discrete_inicon}

On day $0$, we assume that the average arrival flow-rate during interval $m$ ($m=1,\cdots,M$) is given by $f_m^0$. The initial queue length at $t_0$ is $\delta_0^0=0$. Then from \refe{discrete-queue}, the queue length at $t_0+m\dt$ ($m=1,\cdots,M$), $\delta_m^0$, is updated by
\bqn
\delta_m^0&=&\max\{0, \delta_{m-1}^0+ (f_m^0-C)\dt\},
\eqn
and the arrival flow-rate during interval $m$, $g_m^0$, is given by the discrete version of \refe{pq1} or \refe{arrival-rate}:
\bqn
g_m^0&=&f_m^0-\frac{\delta_m^0-\delta_{m-1}^0}\dt=\min\{ \frac{\delta_{m-1}^0}\dt+f_m^0, C\}.
\eqn
We denote the corresponding cumulative departure and arrival flows at $t_0+m \dt$ ($m=0,\cdots,M$) by $F_m^0$ and $G_m^0$, respectively; they are calculated by ($m=0,\cdots,M-1$):
\bqn
F_0^0&=&0, \quad F_{m+1}^0=F_m^0+f_m^0 \dt,\\
G_0^0&=&0, \quad G_{m+1}^0=G_m^0+g_m^0 \dt.
\eqn

We introduce a piecewise linear interpolation of $F_m^0$ for $t\in [t_0+m\dt, t_0+(m+1)\dt]$ ($m=0,\cdots,M-1$) as
\bqn
\tilde F^0(t)&=&F_m^0+(t-t_0-m \dt) f_m^0.
\eqn
Then the queueing time for vehicles arriving at $t_0+m\dt$ ($m=0,\cdots,M$), $\Upsilon_m^0$, is given by
\bqn
\Upsilon_m^0&=&\max \{y| \tilde F^0(t_0+m\dt-y)=G_m^0 \}. \label{queuetime_cal}
\eqn
That is, the cumulative arrival flow at $t_0+m\dt$ equals the cumulative flow at $t_0+m\dt-\Upsilon_m^0$. Equation \refe{queuetime_cal} can be numerically solved in the following steps: (i) If $F_m^0=G_m^0$, $\Upsilon_m^0=0$ from Lemma \ref{lemma:zeroqueue}. (ii) Otherwise, if $F_m^0>G_m^0$, we find the largest $m'$ such that $F_{m'}^0<G_m^0$; then
\bqn
\Upsilon_m^0&=&(m-m')\dt - \frac{G_m^0-F_{m'}^0}{f_{m'}^0}.
\eqn
With $\Upsilon_m^0$, we can calculate the costs at $t_0+m\dt$ and day 0, denoted by $\phi_{1,m}^0$, $\phi_{2,m}$, and $\phi_m^0$, respectively from \refe{def:cost1}, \refe{def:cost2}, and \refe{def:totalcost}.

From \refe{def:rho_int}, we can calculate the average initial imaginary density in cell $i$, $k_i^0$, as
\bqn
k_i^0&=&\frac 1 \dx \int_{y=(i-1)\dx}^{i\dx} k(0,y) dy \nonumber\\
&=&\frac 1 \dx\int_{t_{1,i-1}}^{t_{1,i}} g(0,t) dt+\frac 1 \dx\int_{t_{2,i}}^{t_{2,i-1}} g(0,t) dt \nonumber\\
&=&\frac \dt \dx \sum_{t_{m-\frac 12} \in (t_{1,i-1}, t_{1,i}] \cup [t_{2,i},t_{2,i-1}) } g_m^0 .
\eqn

\subsection{The Cell Transmission Model}\label{section:CTM}

At $j\Delta r$ (or simply $j$th day, if we assume $\Delta r$ is a day), the average density inside cell $i$ between $(i-1)\dx$ and $i\dx$ is denoted by $k_{i}^j$. Then the LWR model, \refe{lwr_sc}, with initial densities $k_i^0$ calculated in the preceding subsection,  can be discretized into the following Cell Transmission Model, which is an extension of the Godunov method \citep{daganzo1995ctm,lebacque1996godunov}: ($i=0,-1,\cdots, -I$)
\ben
\item The demand and supply in cell $i$  are $d_i^j=D(k_i^j)$ and $s_i^j=S(k_i^j)$.
\item From the boundary conditions in \refe{bc_sp}, the downstream supply in cell $1$ is $s_1^j=0$, and the external demand is $d_{-I}^j=0$.
\item The boundary flux at $i\dx$ is given by
\bqn
q_i^j&=&\min\{d_i^j, s_{i+1}^j\}.
\eqn
\item The imaginary density on the next day $(j+1)\Delta r$  in cell $i$ is updated from the following conservation equation:
\bqn
k_i^{j+1}&=&k_i^j+\frac{\Delta r}{\dx} (q_{i-1}^j -q_{i}^j).
\eqn
\een

The Cell Transmission Model is well-defined if and only if $k_i^j\in [0,\kappa]$ for all $i$ and $j$. 
\begin{theorem}
The Cell Transmission Model is well-defined if and only if the following CFL condition is satisfied \citep{courant1928CFL}:
\bqn
\frac {\dx} {\Delta r} \geq &\min\{u,w\}. \label{CFL_condition}
\eqn
\end{theorem}
{\em Proof}. The proof is straightforward and omitted. \eop

\subsection{Calculation of arrival and departure flow-rates}

Once the densities are updated, we can update the corresponding arrival and departure flow-rates.

On day $j$, if $t_{m-\frac 12} \in (t_{1,i-1}, t_{1,i}] \cup [t_{2,i},t_{2,i-1})$, from \refe{density2arrival} we have
\bqn
g_m^j=\frac{\beta \gamma}{\beta+\gamma} k_i^j.
\eqn

We denote $I_*^j$ such that $k_i^j=\kappa$ for $i>-I_*^j$, but $k_i^j<\kappa$ for $i\leq -I_*^j$. That is, cells $0,-1,\cdots, -I_*^j+1$ form the largest jammed scheduling payoff interval containing $0$. 
Then for cell $i\leq -I_*^j$, and time interval $m$, where $t_{m-\frac 12} \in (t_{1,i-1}, t_{1,i}] \cup [t_{2,i},t_{2,i-1})$, we set $f_m^j=g_m^j$.
Further from \refe{equilibrated_departure_flowrate}, we have
\bqn
f_m^j&=&\cas{{ll} \frac{C}{1-\frac \beta \alpha}, & t_{m-\frac 12} \in [t_{1,-I_*^j}, t_*^j);\\
	\frac{C}{1+\frac \gamma \alpha}, & t_{m-\frac 12} \in [t_*^j, t_{2,-I_*^j}),}
\eqn
where $t_*^j=\frac \beta \alpha t_{1,-I_*^j}+(1-\frac\beta\alpha) t_*=-\frac \gamma \alpha t_{2,-I_*^j}+(1+\frac\gamma\alpha) t_*$, and calculate $F_m^j$ and $G_m^j$ from $F_{m+1}^j=F_m^j+f_m^j \dt$ and $G_{m+1}^j=G_m^j+g_m^j \dt$.

The cost functions on day $j$ will be
\bqn
\phi(j\Delta r, t_0+m\dt)&=&\cas{{ll} \phi_2(t_0+m\dt), & t_0+m\dt\leq t_{1,-I_*^j} \m{ or } t_0+m\dt\geq t_{2,-I_*^j};\\
	\phi_2(t_{1,-I_*^j})=\phi_2(t_{2,-I_*^j}), & \m{otherwise;}}\\
\phi_1(j\Delta r, t_0+m\dt)&=&\phi(j\Delta r, t_0+m\dt)-\phi_2(t_0+m\dt).
\eqn

\section{Equilibrium state and its stability}
In this section we study the equilibrium state for the dynamical system of scheduling payoff choice, \refe{lwr_sc}, and its stability property.

\subsection{Equilibrium state}
Given initial and boundary conditions, the fluid will settle down to the bottom part of the single tube in \reff{single_tube}(b) and reach an equilibrium state. Therefore, in the equilibrium state, the imaginary density is given by
\bqn
k^*(x)&=&\cas{{ll} \kappa, & - L^* \leq x\leq 0; \\ 0, &\m{otherwise},}
\eqn
where
\bqn
L^*&=&\frac{N} {\kappa}. \label{def:Lstar}
\eqn
That is, all vehicles have their scheduling costs not greater than $L^*$. The corresponding smallest and largest arrival times are $t_1(L^*)$ and $t_2(L^*)$.
The total cost for all drivers is $L^*$, and the departure flow-rates can be computed from \refe{equilibrated_departure_flowrate}. We refer to the equilibrium state as the scheduling payoff user equilibrium (SPUE), in which the corresponding arrival and departure flow-rates are determined by the splitting and cost balancing procedures in Section \ref{section:costbalancing}.

\begin{theorem}
	With the cost balancing principle, the equilibrium state (SPUE) of the LWR model is equivalent to both ATUE and DTUE.
	\end{theorem}
{\em Proof}. In the SPUE, the first and last arrival times share the same scheduling cost, which is smaller than that at any other time. With the cost balancing principle we can see that all vehicles arriving between the two times have the same cost. Thus this is also the ATUE. 

When a system is in ATUE, no vehicles can find an under-utilized arrival time with a smaller scheduling cost. Thus the LWR model has to be in the SPUE; otherwise, there can be positive flows.

Therefore SPUE together with the cost balancing principle, the SPUE is equivalent to ATUE. Further from Theorem \ref{theorem:DTUE-ATUE-equivalence}, these three equilibria are equivalent.
 \eop

\subsection{Stability}\label{section:stability}
We set the cell size to be $L^*$ and consider a small perturbation to the equilibrium state as
\bqn
k^*(x)+\epsilon_0(r,x)&=&\cas{{ll} \kappa-\epsilon(r), & - L^* \leq x\leq 0;\\ \epsilon(r), & - 2 L^* \leq x\leq - L^*; \\ 0, &\m{otherwise}.}
\eqn
Thus the density is always zero in cells other than 0 and -1. 

We denote the average densities at $r$ in cells 0 and -1 by $k_0(r)$ and $k_{-1}(r)$. Due to the conservation law, we have
\bqn
k_0(r)+k_{-1}(r)&=& \kappa.
\eqn
Letting $k_{-1}(r)=\epsilon(r)$, we have $k_0(r)=\kappa-\epsilon(r)$ and $\epsilon(0)=\epsilon$.

Assuming that $\epsilon(r)$ is always small, we can calculate the demands and supplies in cells 0 and -1 as
\bqs
d_0(r)&=&u \kappa_c,\\
s_0(r)&=&w \epsilon(r),\\
d_{-1}(r)&=&u \epsilon(r),\\
s_{-1}(r)&=&u \kappa_c.
\eqs
In addition $d_{-2}(r)=0$ and $s_{1}(r)=0$. Thus the boundary fluxes are
\bqs
q_0(r)&=&\min\{d_0(r),s_1(r)\}=0,\\
q_{-1}(r)&=&\min\{d_{-1}(r),s_0(r) \}=\min\{u,w\}  \epsilon(r),\\
q_{-2}(r)&=&\min\{d_{-2}(r),s_{-1}(r) \}=0.
\eqs
Then the LWR model, \refe{lwr_sc}, can be approximated by the following ordinary differential equations in cells 0 and -1:
\bqs
\der{}{r} k_0(r)&=& \frac 1{L^*} (q_{-1}(r)-q_0(r) ),\\
\der{}{r} k_{-1}(r)&=& \frac 1{L^*} (q_{-2}(r)-q_{-1}(r) ),
\eqs
which both lead to
\bqn
\der{}{r} \epsilon(r) &=&- \frac 1{L^*} \min\{u,w\}  \epsilon(r). \label{epsilon-system}
\eqn

\begin{theorem}
	The ordinary differential equation, \refe{epsilon-system}, is asymptotically stable and converges to $\epsilon(\infty)=0$.
\end{theorem}
{\em Proof}. \refe{epsilon-system} is asymptotically stable since the eigenvalue of the linear system is $- \frac 1{L^*} \min\{u,w\}$, which is negative. In addition, it is
solved by
\bqs
 \epsilon(r)&=& \epsilon e^{- \frac 1{L^*} \min\{u,w\} r },
\eqs
which converges to $0$ for $r\to \infty$.
\eop

Since \refe{epsilon-system} approximates the LWR model, \refe{lwr_sc}, we expect that the latter is also stable and converges to the SPUE. That is, the day-to-day dynamics of scheduling payoff, arrival time, and departure time choices are asymptotically stable.

\section{Numerical example} \label{section:numerical}

In this section, we set $N=3600$ veh, $C=1800$ vph, $\alpha=50$ \$/hr, $\beta=25$ \$/hr, and $\gamma=100$ \$/hr. Thus the jam density for the single tube is $\kappa=90$ veh/\$ from \refe{def:kappa},  and the final equilibrium state's length is  $L^*=40$ \$ from \refe{def:Lstar}. 

The ideal arrival time $t_*=0$.  
The study time period is set as $[-4, 1]$ hr, where $\phi_2(-4)=\phi_2(1)=100$ \$. Thus the length of the single tube is $L=100$ \$, and $x\in [-L,0]$ in the LWR model. We let $u=w=1$ \$/day. Thus $\kappa_c=45$ veh/\$. 

We set $\dt=\frac 1{1000}$ hr, $\dx=\frac 12$ \$, and $\Delta r= \frac 12$ day, which satisfy the CFL condition \refe{CFL_condition}. Thus the study time period is divided into 5000 intervals, and the imaginary road is split into 200 cells. Corresponding to the 200 cells, the 400 time intervals have the following boundaries: $t_{1,i}=0.04 i \dx$ and $t_{2,i}=-0.01 i\dx$ for $i=-200,\cdots, 0$.

\bfg\bc
\includegraphics[width=6in]{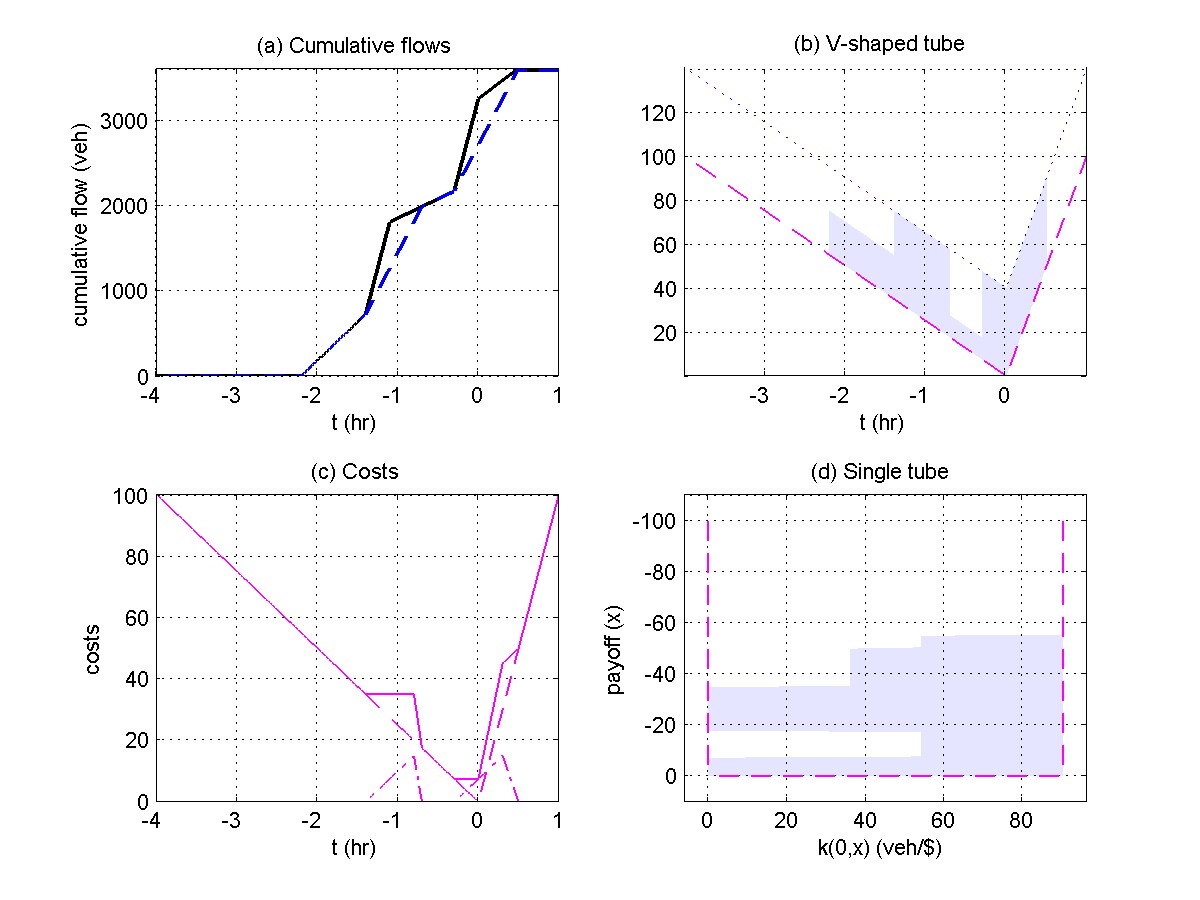}\caption{The initial conditions for the numerical example: (a) Cumulative flows $F(t)$ (solid line) and $G(t)$ (dashed line); (b) The arrival flow-rates at different times in a V-shaped tube; (c) Costs $\phi_1(t)$ (dash-dotted line), $\phi_2(t)$ (dashed line), and $\phi(t)$ (solid line); (d) Imaginary densities at different scheduling payoffs in a single tube (or imaginary road). }\label{numerical_initial}
\ec\efg

The initial departure flow-rate is given by
\bqs
f(0,t)&=&\cas{{ll} \frac 12 C, & t\in(-2.2, -1.4];\\ 2 C, & t\in (-1.4, -1.1]; \\\frac 14 C, & t\in (-1.1, -0.3]; \\
	2C, &t\in (-0.3, 0]; \\ 0.4 C, & t\in (0,0.5]; \\
	0, &\m{otherwise},
 }                    
\eqs
which is similar to the departure flow-rate shown in \reff{bottleneck_costs}.
With the numerical methods in Section \ref{section:discrete_inicon}, we can calculate the arrival flow-rate and cumulative flows, the corresponding costs, and the initial imaginary density at different scheduling payoffs. The initial conditions are shown in \reff{numerical_initial}. 

\bfg\bc
\includegraphics[width=5in]{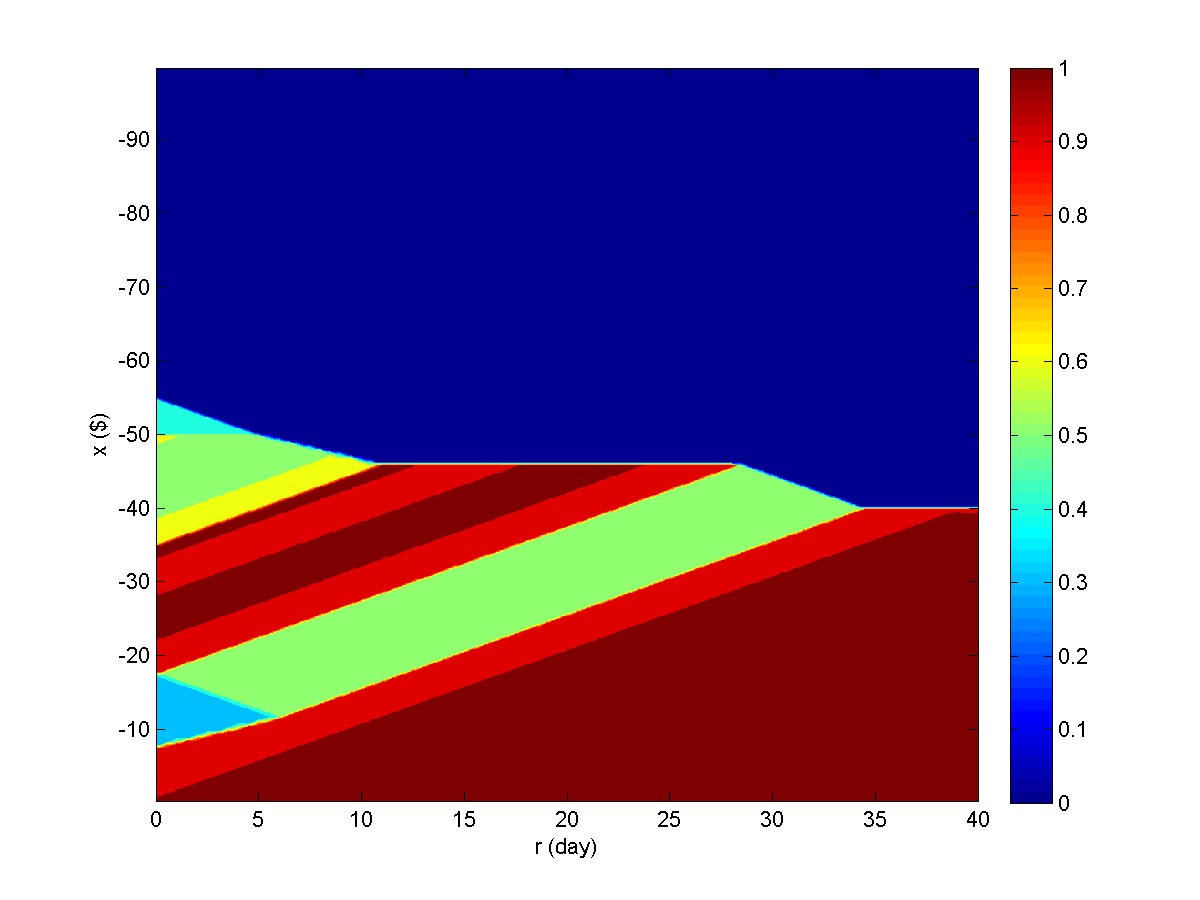}
\caption{Evolution of normalized imaginary density}\label{fig:density_evolution}
\ec\efg

Solving the Cell Transmission Model in Section \ref{section:CTM}, we obtain the  flow density normalized with respect to the jam density, $\kappa$, as shown in \reff{fig:density_evolution}. From the figure we can see that the flow indeed moves toward the bottom part of the single tube at $x=0$ until it reaches the equilibrium state on day 40, when no further movement is possible. During the process, backward traveling waves form inside the tube with the speed of $-w$. The convergence of the numerical result confirms the stability of the LWR model, \refe{lwr_sc}, as demonstrated in Section \ref{section:stability}.

\bfg\bc
\includegraphics[width=6in]{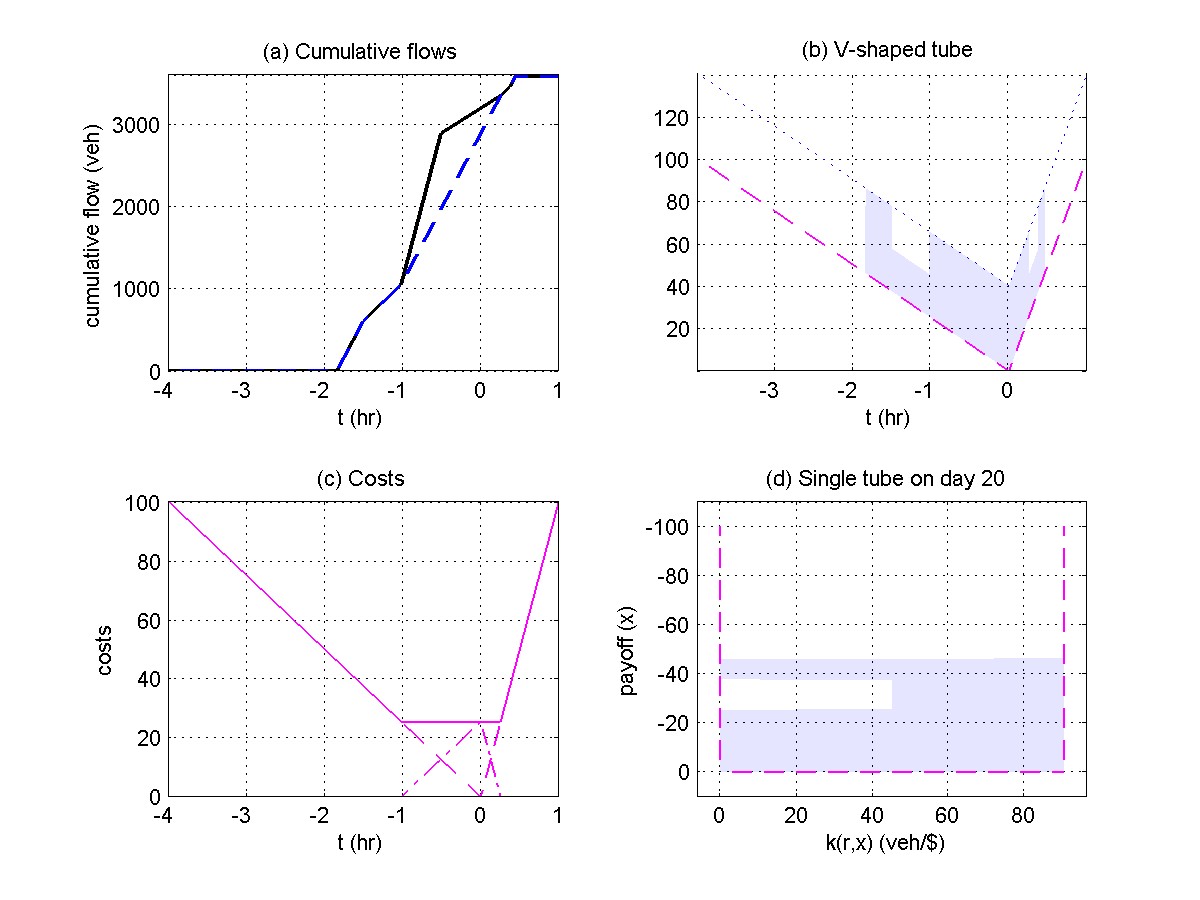}\caption{Numerical solutions on day 20: (a) Cumulative flows $F(t)$ (solid line) and $G(t)$ (dashed line); (b) The arrival flow-rates at different times in a V-shaped tube; (c) Costs $\phi_1(t)$ (dash-dotted line), $\phi_2(t)$ (dashed line), and $\phi(t)$ (solid line); (d) Imaginary densities at different scheduling payoffs in a single tube. }\label{numerical_day20}
\ec\efg

\reff{numerical_day20} and \reff{numerical_day40} respectively illustrate the solutions on days 20 and 40. Compared with the initial conditions, the solution on day 20 is closer to the SPUE state. In particular, the solutions on day 40 verify the existence and stability of the SPUE, ATUE, and DTUE.

\bfg\bc
\includegraphics[width=6in]{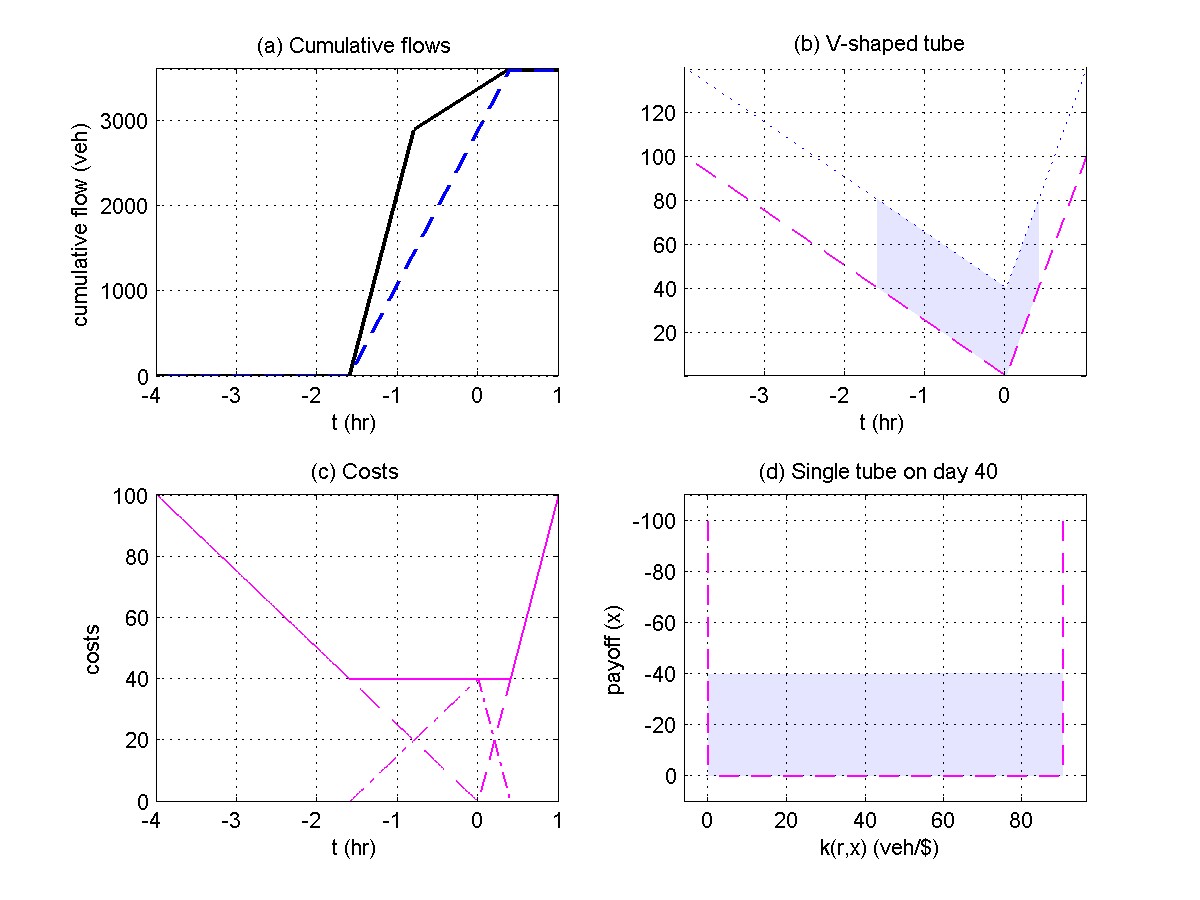}\caption{Numerical solutions on day 40: (a) Cumulative flows $F(t)$ (solid line) and $G(t)$ (dashed line); (b) The arrival flow-rates at different times in a V-shaped tube; (c) Costs $\phi_1(t)$ (dash-dotted line), $\phi_2(t)$ (dashed line), and $\phi(t)$ (solid line); (d) Imaginary densities at different scheduling payoffs in a single tube. }\label{numerical_day40}
\ec\efg

\section{Conclusion}
In this paper we presented a new day-to-day dynamical system model for drivers' departure time choice at a single bottleneck. We first defined the within-day dynamics with a point queue, associated costs, and the departure time and arrival time user equilibria. Based on three behavioral principles, we showed that that drivers choose their scheduling payoffs first, then arrival times, and finally departure times. We then presented conceptual models for both arrival time and scheduling payoff choices; in particular, with a single tube model, we showed that the nonlocal arrival time shifting problem is equivalent to the local scheduling payoff shifting problem. After introducing the imaginary density variable, we applied the Lighthill-Whitham-Richards model to describe the day-to-day dynamics for scheduling payoff choice and presented the splitting and cost balancing procedures for arrival and departure time choices. We also presented the corresponding discrete models for their numerical solutions. We theoretically proved that the equilibrium state (scheduling payoff user equilibrium) of the LWR model corresponds to the user equilibrium and is stable. We used one numerical example to demonstrate the effectiveness and stability of the new day-to-day dynamical system model.

Such a stable dynamical system can be used as a computational method to calculate the arrival and departure time choices, to understand drivers' learning processes in deciding arrival and departure times, to study the impacts of information provision, and to develop effective and efficient pricing and other control methods.

This study is the first step for understanding stable day-to-day dynamics for departure time choice. 
Even though plausible, drivers' choice behaviors as well as the corresponding conceptual and mathematical models in this study are subject to further theoretical and empirical analyses. For examples, are the LWR model, \refe{lwr_sc}, and the triangular fundamental diagram, \refe{lwr_tri_fd}, reasonable for people's scheduling payoff choice? Is there a better model than \refe{density2arrival} for drivers' arrival time choice? How accurate is \refe{equilibrated_departure_flowrate} for drivers' departure time choice?

Mathematically, the proof of the day-to-day dynamical system model's stability in Section \ref{section:stability} is based on a locally approximate ordinary differential equation, which can be viewed as Lyapunov's first method. Even though the stability has been verified numerically in Section \ref{section:numerical}, we will be interested in directly proving the stability of the LWR model, which is a partial differential equation subject to boundary conditions. In particular, we will be interested in using the Lyapunov's second method to prove the stability.

Some extensions are possible with the new model. The day-to-day dynamical system can be extended for the dynamic traffic assignment problem with simultaneous departure time and route choices \citep{Friesz1993due}. For example, the V-shaped tube model of arrival time choice illustrated in \reff{bottleneck_tube} can be extended for more complicated scheduling cost functions, such as a W-shaped function. The dynamical system model can also be extended for multi-class traffic systems, where drivers have different values of time, stochastic distribution of ideal arrival times, and/or scheduling cost functions. We will also be interested in developing a day-to-day dynamical system model of departure time choice for dynamic/varying capacities, parallel roads, and general networks. We will examine how different number of commuters can impact the system's performance. The model could lead to an effective method for computing the elusive user equilibrium under complex conditions \citep{smith1984existence,daganzo1985uniqueness}.
Finally, with a better understanding drivers' behaviors at such a bottleneck, we can provide more efficient guidance for human drivers and connected and autonomous vehicles regarding departure time choice at various bottlenecks.

\section*{Acknowledgments}
Discussions with Prof. Hai Yang of the Hong Kong University of Science and Technology, Emeritus Professor Pete Fielding and Ms. Irene Martinez of UC Irvine, Prof. Kentaro Wada of the University of Tokyo, and Prof. Takamasa Iryo of Kobe University are acknowledged. The views and results are the author's alone.

\pdfbookmark[1]{References}{references}

\end {document}